\RequirePackage[leqno]{amsmath}
\documentclass[leqno]{amsart}
\makeatletter
\usepackage[utf8]{inputenc}     
\usepackage[T1]{fontenc}
\usepackage{multicol}
\usepackage{enumitem}
\usepackage{color, graphics}
\usepackage[demo]{graphicx}
\usepackage{efbox,graphicx}
\efboxsetup{linecolor=gray,linewidth=0.8pt}
\usepackage{floatrow}
\usepackage{pgfplots}
\usepackage{graphicx,multicol,latexsym,amsmath,amssymb}
\usepackage{pgfplots}
\usepackage{tikz} 
\usepackage{amsfonts}
\usepackage{enumitem}
\setlist[itemize]{topsep=0pt,after=\vspace{1.5\baselineskip}}
 \usepackage[colorlinks=true]{hyperref}
\usepackage{empheq}
\usepackage{tabularx}
\usepackage[leftcaption]{sidecap}
\sidecaptionvpos{figure}{m}
\usepackage{tikz}
\usepackage[margin=1in]{geometry}
\usepackage{float}
\usepackage{subcaption}
\usetikzlibrary{backgrounds,automata}

\newcount\rc@count
\let\rc@clearconstantlist\empty
\newcommand\rc@clearconstant[1]{\global\expandafter\let\csname rc@const@#1\endcsname\undefined}
\newcommand\resetconstants[1]{%
    \def\rc@constname{#1}
    \global\rc@count=1\relax 
    \bgroup 
        \let\\\rc@clearconstant 
        \rc@clearconstantlist
        \global\let\rc@clearconstantlist\empty 
    \egroup
}
\newcommand\const[1]{%
    \@ifundefined{rc@const@#1}{%
        \expandafter\xdef\csname rc@const@#1\endcsname{%
           \noexpand\rc@useconst{\rc@constname}{\the\rc@count}%
        }%
        \g@addto@macro\rc@clearconstantlist{\\{\mathrm{#1}}}%
        \global\advance\rc@count1\relax
    }{}%
    \csname rc@const@#1\endcsname
}

\newcommand\rc@useconst[2]{{#1}\textsubscript{#2}}
\setlist[itemize]{noitemsep, topsep=0pt}

\def\R{\mathbb R} \def\N{\mathbb N}

\def\R{\mathbb R} \def\N{\mathbb N} 
\def\TM{T_{\rm{max}}} 
\def
\@cite
#1#2{[#1\if@tempswa, #2\fi]}
\makeatother

\addcontentsline{toc}{section}{Abstract}
\newtheorem{theorem}{Theorem}[section]
\newtheorem{corollary}[theorem]{Corollary}

\newtheorem{lemma}[theorem]{Lemma}

\newtheorem{remark}{Remark}

\newcounter{cnstcnt}

\title[Boundedness in chemotaxis models with local, nonlocal, and gradient effects]{
Chemotaxis models with mixed mechanisms: boundedness in growth-dominated regimes}
\author[Tongxing Li, Silvia Frassu and Giuseppe Viglialoro]{}
\keywords{Chemotaxis, Global existence, Nonlocal terms, Gradient terms, Boundedness. \\
\textit{$^\star$Corresponding author}: silvia.frassu@unica.it}

\begin{document}

\maketitle
\centerline{\scshape{\scshape{Tongxing Li$^{\natural}$, Silvia Frassu$^{\sharp,\star}$ \and Giuseppe Viglialoro$^{\sharp}$}}}
\medskip
{
\medskip
\centerline{$^{\natural}$School of Control Science and Engineering} 
\centerline{Shandong University}
\centerline{Jinan, Shandong, 250061 (P. R. China)}
\medskip
}
{
\medskip
\centerline{$^\sharp$Dipartimento di Matematica e Informatica}
\centerline{Universit\`{a} degli Studi di Cagliari}
\centerline{Via Ospedale 72, 09124. Cagliari (Italy)}
\medskip
}
\bigskip
\begin{abstract}
We investigate a class of zero-flux chemotaxis-growth systems featuring nonlinear local and nonlocal reaction terms, as well as a gradient-dependent damping, given by
\begin{equation*}
\begin{cases}
u_t = \Delta u - \chi \nabla \cdot (u \nabla v) + a u^\rho - b \left( \int_\Omega u^\beta \right)^\delta - c |\nabla u|^\gamma, & \text{in } \Omega \times (0,\TM), \\
\tau v_t = \Delta v - v + u, & \text{in } \Omega \times (0,\TM), \\
\end{cases}
\end{equation*}
where $\chi, a, b, c > 0$, $\rho, \beta, \delta, \gamma \geq 1$, and $\tau \in \{0,1\}$. Assuming that $\Omega \subset \mathbb{R}^n$ ($n\geq 1$) is a bounded domain with smooth boundary and that the initial data $(u_0, \tau v_0)$ are sufficiently regular, we prove the existence of global (i.e., $\TM=\infty$), uniformly bounded classical solutions under suitable structural conditions on the exponents defining the reaction term and the dimension $n$.

A central insight of our analysis is that, in contrast to classical scenarios where $\rho \leq \beta$ typically ensures automatic control of the total mass $\int_\Omega u$, the regime $\rho > \beta$ makes that  the structure of the source term alone does not suffice to guarantee mass boundedness. Accordingly, the global boundedness of solutions is achieved via a two-tiered strategy: first, we identify a parameter regime in which the total mass remains uniformly bounded in time; second, under further constraints on the model parameters, we establish uniform-in-time boundedness of classical solutions in stronger norms, specifically in $L^\infty(\Omega)$.

These results offer further insight into the interplay between chemotactic aggregation, gradient-driven dissipation, and nonlocal reaction effects, contributing to the analysis of blow-up prevention in structured chemotaxis models.
\end{abstract}
\resetconstants{c}
\section{Introduction}\label{Intro}
\subsection{Reaction-Diffusion-Transport phenomena: the Keller--Segel model}
Nonlinear partial differential equations are widely used to model complex biological and chemical processes, including population growth, migration, and chemical signaling. These systems often involve local reactions, spatial diffusion, transport influenced by gradients, and global feedback mechanisms. Nonlocal effects, such as total population or concentration-dependent regulation, play a crucial role in capturing collective behaviors like quorum sensing or resource limitation. The class of equations studied here encompasses such interactions and falls within the broader framework of reaction-diffusion-transport models.

A reaction-diffusion-transport phenomenon is generally formulated by a partial differential equation of the type:
\begin{equation}\label{TaxisEquationTransport}
u_t = \Delta u - \nabla \cdot (u w) + R(x,t,u,\nabla u).
\end{equation}
Above, 
\begin{itemize}
\item [$\circ$] \( u = u(x,t) \) is the unknown, as for instance a population density concentration in a domain at a position $x$ and at the time $t$;
  \item [$\circ$] \( \Delta u \) represents the diffusion term, accounting for random spatial motion;
  \item [$\circ$] \( - \nabla \cdot (u w) \) is the transport (or advection) term, which describes how the quantity $u$ is carried along by the vector flow $w=w(x,t)$; 
  \item [$\circ$] \( R(x,t,u,\nabla u) \) is the reaction term, possibly depending on \( u \), $\nabla u$, \( x \) and \( t \), models growth/decay mechanisms.
\end{itemize}
In the classical Keller--Segel model, the vector flow is the gradient of a further unknown $v=v(x,t)$ coupled to the equation of $u$ itself. In the specific, the system
\begin{equation}\label{KS}
 u_t = \Delta u - \nabla\cdot(u\nabla v), \qquad v_t = \Delta v - v + u
\end{equation}
was introduced to describe aggregation in starving slime moulds \cite{K-S-1970,BellomoEtAl}, where individual cells (density $u$) bias their random movement toward higher concentrations of a self-produced chemical signal ($v$).
\subsection{Analysis of the Keller--Segel system}
From a mathematical perspective, the Keller--Segel system \eqref{KS} is well-known for exhibiting both global boundedness and finite-time blow-up of solutions, depending crucially on the spatial dimension, the initial mass, and the structure of initial data. For $n=1$, all solutions are globally bounded. In contrast, for $n \geq 3$ and any prescribed $m>0$, there exist initial data with mass $\int_\Omega u_0 = m$ that lead to blow-up in finite time. The critical case arises in two dimensions: if $m < 4\pi$, diffusion dominates and solutions remain global; if $m > 4\pi$, blow-up can occur. For a comprehensive overview, we refer to the surveys \cite{BellomoEtAl,HerreroVelazquez,Nagai,WinklAggre}.
\subsection{The Keller--Segel system with logistic source terms of different types}
The inclusion of logistic-type source terms in \eqref{KS} introduces superlinear damping mechanisms that are expected to counteract aggregation and prevent blow-up. Specifically, terms of the form \( f(u) = au - bu^\beta \), with \( a, b > 0 \) and \( \beta > 1 \), model population growth limited by density-dependent mortality, corresponding to the choice \( R(x,t,u,\nabla u) = au - bu^\beta \) in the general framework \eqref{TaxisEquationTransport} \cite{verhulst, AiTsEdYaMi, MIMURATsujikawa}. For instance, in the case $\beta = 2$, global existence and boundedness have been proven under suitable conditions \cite{TelloWinkParEl, W0}. On the other hand, finite-time blow-up persists for certain parameter regimes: blow-up has been demonstrated for $\beta=\beta(n) > 1$ (and close to $1$) when $v$ satisfies certain elliptic equations in \cite{WinDespiteLogistic,Winkler_ZAMP-FiniteTimeLowDimension}, with refinements in \cite{FuestCriticalNoDEA}.

Although the logistic source term above models a reaction that depends solely on population density, in reality, population growth or decay may also be influenced by additional factors intrinsically linked to the density itself. At the macroscopic level, behavior effects of species translate into additional loss terms depending on $\nabla u$ in the reaction $R$ of \eqref{TaxisEquationTransport}, yielding for $\chi, c>0$:
\begin{equation}\label{KS-Gradient} 
u_t = \Delta u - \chi \nabla\cdot(u\nabla v) + au^\rho - bu^\beta - c|\nabla u|^\gamma, \quad \textrm{ with } \beta>\rho\geq 1, \gamma\geq 1. 
\end{equation}
In the same spirit, and with an analogous interpretation of the growth and decay mechanisms discussed above, models involving reaction terms with nonlocal effects as
\begin{equation}\label{KS-NonLocal} 
u_t = \Delta u - \chi \nabla \cdot (u \nabla v) + au^\rho - bu^\rho \int_\Omega u^\beta,\quad \textrm{ with } \rho,\beta\geq 1,
\end{equation}
or also others combining nonlocal and gradient nonlinearities, so leading to
\begin{equation}\label{KS-NonLocalGradient}
u_t = \Delta u - \chi \nabla \cdot (u \nabla v)+au^\rho - b u^\rho \int_\Omega u^\beta - c|\nabla u|^\gamma, \quad \textrm{ with } \rho,\beta,\gamma \geq 1,
\end{equation}
have been lately using. 

In this context, it is not our aim to provide a detailed analysis of every existing result concerning chemotaxis models involving cell equations of the forms \eqref{KS-Gradient}, \eqref{KS-NonLocal} and \eqref{KS-NonLocalGradient}. Doing so would shift the focus away from our primary objective. Nevertheless, we provide a set of references where the reader can find comprehensive analyses regarding boundedness, blow-up, and the asymptotic behavior of solutions to models incorporating the aforementioned formulations. See \cite{IshidaLankeitVigliloro-Gradient,LiEtAlStudApplMath2025,BianEtAlNonlocal,BianChenEvangelosLinearNonlocal,
TaoFang-NonlinearNonlocal,Tello2021,DiazFuentesETAl-AAM-2025}. (As a real-world application of the previously discussed models with those sources, taxis systems can describe phenomena such as predator movement toward prey concentrations \cite{kareiva_odell}. Even in the absence of explicit interspecies interaction terms, similar formulations arise when incorporating spatial memory effects \cite{fagan2013memory}, or when modeling attraction to previously visited locations through weak memory kernels in random walk frameworks \cite{shi_shi_wu_spatial_memory}.)
\subsubsection{The Reaction-Diffusion-Transport phenomena as the natural extension of Reaction-Diffusion mechanisms}
As previously discussed, it is well known that purely diffusive models admit global-in-time solutions. In contrast, the introduction of a convective (transport) term can lead to blow-up phenomena and instabilities. In particular, the transport-reaction-diffusion models considered above represent a natural extension of models involving only diffusion and reaction.

In this context, the influence of gradient-dependent and nonlocal effects on blow-up phenomena has been extensively investigated in nonlinear parabolic equations. Such models naturally emerge in biological settings, for example, in the dynamics of a single species exhibiting population-dependent behavioral responses or facing resource constraints influenced by the overall density. Mathematically, this leads to equations of the type
\begin{equation*}
u_t = \Delta u +a u^\rho- c|\nabla u|^\gamma, \quad \textrm{or for proper $f$}\quad u_t = \Delta u + a u^\rho + f(u) \left( \int_\Omega u^\beta \right)^\delta,
\end{equation*}
and classical works dealing with properties of related solutions are, for instance,  \cite{chipot_weissler,kawohl_peletier,fila,SoupBio,PolacikQuittnerSouplet2007,HuYin1995}. In particular, in these studies threshold relations between the exponents $\rho,\beta,\delta$ that separate global existence from finite-time blow-up may be found in some specific scenarios; see also \cite[$\S$IV and $\S$V]{quittner2007superlinear}.
%
%

\section{Main motivations: on the evolution of $\int_\Omega u(x,t)dx$ over time}\label{MotivationSection}
\subsection{The natural uniform-in-time boundedness of the mass}
Although not explicitly stated, the mathematical analysis of the equations describing the aforementioned phenomena, with or without transport, is naturally carried out by considering them on a space-time domain \(\Omega \times (0, T_{\max})\), where \(\Omega\) is a bounded and smooth domain in \(\mathbb{R}^n\), and \(T_{\max} \in (0, \infty]\) denotes the maximal existence time of solutions. This setting is completed by prescribing appropriate boundary conditions on \(\partial \Omega \times (0, T_{\max})\) and initial data \(u(x,0) = u_0(x)\) in \(\Omega\).

With regard to boundary restrictions, homogeneous Neumann conditions are considered the most natural within the framework of chemotactic systems. In this context, the presented results hold for initial boundary value problems subject to $u_\nu = v_\nu = 0$ on $\partial \Omega \times (0,\TM)$, where $\nu$ denotes the outward normal vector to the boundary of $\Omega$, $\partial \Omega$.

By analyzing equations \eqref{KS-Gradient}, \eqref{KS-NonLocal} and \eqref{KS-NonLocalGradient}, and in particular by considering the interplay between the nonlinear terms and the relative magnitudes of the exponents involved, one can deduce that any nonnegative solution \( u=u(x,t) \), defined on the space-time domain \( \Omega \times (0, \TM) \), satisfies a uniform-in-time bound on its total mass. More precisely, the \( L^1(\Omega) \)-norm of \( u \), given by $
\int_\Omega u$, remains uniformly bounded for all \( t \in (0, \TM) \), and having such a property at disposal is essentially necessary in order to carry out subsequent analyses, concerning, for example, the uniform boundedness of the function $u$ itself.

To show that the $L^1(\Omega)$-norm of $u$ is uniformly bounded for the equation \eqref{KS-Gradient}, one can proceed in various ways; however, here we prefer to provide a constructive geometric justification. (We control the reaction term that is eventually concave by means of a tangent line at an appropriate point.) In particular we will show that for all $a,b>0$ and $\beta>\rho\geq 1$ there are $C_0>0$ and $C_1>0$ such that
\begin{equation}\label{Eq1}
\psi(s):=a s^\rho - b s^\beta \leq C_0 - C_1 s \quad \text{for all } s > 0.
\end{equation}
As a consequence, by integrating \eqref{KS-Gradient} over $\Omega$, the divergence theorem gives, in view of the fact that the diffusion and cross-diffusion terms vanish due to the Neumann boundary conditions and \eqref{Eq1} itself,
\begin{equation*}\label{Eq0}
\frac{d}{dt} \int_\Omega u  = a \int_\Omega u^\rho  - b \int_\Omega u^\beta  - c \int_\Omega |\nabla u|^\gamma \leq a \int_\Omega u^\rho  - b \int_\Omega u^\beta\leq \int_\Omega C_0-C_1\int_\Omega u\quad \textrm{on} \;(0,\TM), 
\end{equation*}
%
leading to $$\int_\Omega u\leq \max\left\{\int_\Omega u_0(x)dx,\frac{C_0|\Omega|}{C_1}\right\} \quad \textrm{on}\quad (0,\TM).$$
Now, since \( \beta > \rho\geq 1 \),  for  
$$c_m=\frac{1}{2} \left(\left(\frac{a}{b}\right)^{\frac{1}{\beta-\rho}}+\left(\frac{a \rho}{b \beta}\right)^{\frac{1}{\beta-\rho}}\right), \quad c_t=\psi(c_m),$$
it is not difficult to check that relation \eqref{Eq1} holds with 
\[C_0=c_t-\psi'(c_m)c_m\quad \textrm{and}\quad C_1=-\psi'(c_m).\]
As to the related analysis for equations  \eqref{KS-NonLocal} and \eqref{KS-NonLocalGradient}, similarly we have by Hölder's inequality
\begin{equation*}
\frac{d}{dt} \int_\Omega u  \leq \int_\Omega u^\rho\left(a-b \int_\Omega u^\beta\right)\leq \int_\Omega u^\rho\left(a-b |\Omega|^{1 - \beta} \left(\int_\Omega u\right)^\beta\right) \quad \textrm{for all } t\in (0,\TM),
\end{equation*}
and we obtain the uniform-in-time boundedness of $\int_\Omega u$ on $(0,\TM)$ thanks to \cite[Lemma 3.3]{ChiyoEtAl2024}.
\subsection{Growth mechanisms prevailing over decay: lack of automatic mass control}
It has been observed that reaction-diffusion systems featuring well-defined logistic-type nonlinearities often exhibit the property that the total mass of the solution (defined as the spatial integral of the unknown function) can be effectively controlled over  $(0,\TM)$.

In the present work we are concerned with a different scenario, in which the total mass is not automatically controlled by the structure of the equation. Rather, it becomes necessary to impose explicit conditions to ensure its boundedness. This situation arises particularly in models where the reaction term includes a superlinear growth component, while the dissipative (or decay) effects, possibly of nonlocal type, are not dominant in general. In particular,  we will consider the following equation
\begin{equation}\label{EqNoMass}
u_t=\Delta u -\chi \nabla \cdot (u \nabla v)+a u^\rho  -b \left(\int_\Omega u^\beta\right)^\delta- c |\nabla u|^\gamma\quad \textrm{on } \Omega \times (0,\TM),
\end{equation}
to be naturally coupled with an equation for $v$ (specified later on).

Crucially, if in equation \eqref{EqNoMass} we assume that $\rho \leq \beta$, arguments analogous to those presented above ensure that the mass cannot grow definitively, so yielding a bound for $\int_\Omega u$ on the interval $(0, \TM)$. In contrast, when
\begin{equation}\label{CondrhoMAGbeta}
\rho > \beta\geq 1,
\end{equation}
this integrability property of $u$ is no longer preserved. In our setting, henceforth, the decay mechanism is not necessarily stronger than the growth at the outset. Therefore, suitable conditions on the coefficients and exponents in the decay term to recover mass control and prevent unphysical growth of the solution must be placed.

Building on the preceding analysis and motivated by the associated discussion, the following section introduces our chemotaxis model, which is rooted in equation \eqref{EqNoMass}, subject to the constraint \eqref{CondrhoMAGbeta}.
\section{The model and main results }
Let us turn our attention to this chemotaxis-growth model
\begin{equation}\label{problem}
\begin{cases}
u_t=\Delta u -\chi \nabla \cdot (u \nabla v)+a u^\rho  -b \left(\int_\Omega u^\beta\right)^\delta- c |\nabla u|^\gamma &{\rm in}\ \Omega \times (0, \TM),\\
\tau v_t=\Delta v-v+u &{\rm in}\ \Omega \times (0, \TM),\\
u_\nu=v_\nu=0 &{\rm on}\ \partial\Omega \times (0, \TM),\\
u(x, 0)=u_0(x), \tau v(x,0)= \tau v_0(x) &x \in \bar{\Omega},
\end{cases} 
\end{equation}
where $\chi, a, b, c>0$, $\beta, \delta, \gamma, \rho \ge 1$, $\tau \in \{0,1\}$.
In addition, we impose the following regularity assumptions:
\begin{equation}\label{reglocal}
\begin{cases}
\textrm{ For some } \alpha\in(0,1) \textrm{ and } n\in \N, \Omega \subset \R^n \textrm{ is a bounded domain of class } C^{2+\alpha}, \textrm{ with boundary } \partial \Omega,\\
u_0, \tau v_0: \bar{\Omega}  \rightarrow \R^+ , \textrm{ with } u_0,  \tau v_0 \in C_\nu^{2+\alpha}(\bar\Omega)=\{\psi \in C^{2+\alpha}(\bar{\Omega}): \psi_\nu=0 \textrm{ on }\partial \Omega\}. 
 \end{cases}
\end{equation}
Under these settings, we are now in a position to state our main result.
\begin{theorem}\label{MainTheorem}
For $\tau\in\{0,1\}$, let $\Omega, u_0, \tau v_0$ comply with hypotheses \eqref{reglocal}. Additionally, let $\chi, a, b, c>0$. Then, whenever 
\begin{equation}\label{MT1}
1 \leq \beta < \frac{2(n+1)}{n}, \; \max\{2, \beta\} < \rho < \frac{2(n+1)}{n}, \; \frac{n\rho}{n+1}<\gamma \leq 2 \quad \textrm{and} 
\quad \delta>\max\left\{\frac{n\gamma-n\rho+\rho\gamma}{n\gamma-n\rho+\beta\gamma}, \frac{\rho}{\beta}\right\},    
\end{equation} 
or
\begin{equation}\label{MT2}
1 \leq \beta < \rho \leq 2, \quad \frac{2n}{n+1}<\gamma \leq 2 \quad \textrm{and} 
\quad \delta>\max\left\{\frac{n\gamma-n\rho+\rho\gamma}{n\gamma-n\rho+\beta\gamma}, \frac{\rho}{\beta}\right\},    
\end{equation} 
problem \eqref{problem} admits a unique classical solution $(u,v)=(u(x,t),v(x,t))$, nonnegative and global, in the sense that $u, v \geq 0 $ for all $(x,t) \in \bar{\Omega}\times (0,\infty)$, and uniformly bounded in time, i.e., $u,v \in L^{\infty}(\Omega \times (0,\infty))$.
\end{theorem} 
Moreover, some alternative conditions ensure boundedness for specific values of $\gamma$ and $\delta$. We make it precise in the following 
\begin{corollary}\label{MainCor}
Under the same hypotheses of Theorem \ref{MainTheorem}, let $\delta=\frac{\rho}{\beta}$ and $\gamma=\rho$. Then there is $C_{GN}>0$ such that the same statement holds whenever 
\begin{equation}\label{MC}
\max\left\{\beta, \frac{2n}{n+1}\right\} < \rho \leq 2 \quad \textrm{and} \quad  b \geq \frac{a 2^{\rho-1}C_{GN}^\rho}{|\Omega|} \left[\frac{\beta \rho}{n\rho+\beta \rho-n \beta} \left(\frac{a 2^{\rho-1}C_{GN}^\rho n (\rho-\beta)}{c(n\rho+\beta \rho-n \beta)}\right)^{\frac{n(\rho-\beta)}{\beta \rho}}+1\right].
\end{equation}
\end{corollary}

\begin{remark}
The results presented in Theorem~\ref{MainTheorem} and Corollary~\ref{MainCor} highlight a crucial structural feature of system~\eqref{problem}: in regimes where the nonlinear growth exponent $\rho$ exceeds decay effects, the emergence of global, bounded solutions is contingent upon the presence of a sufficiently strong damping mechanism.  This damping effect is realized both through the nonlocal absorption term, characterized by the exponent $\delta$, and the nonlinear diffusion of gradient-type, governed by $\gamma$. In particular, the threshold conditions imposed on $\delta$ in \eqref{MT1} and \eqref{MT2}, as well as the sharp balancing condition $\delta = \rho/\beta$ in Corollary~\ref{MainCor}, demonstrate that when the system allows for rapid, possibly superlinear, growth (i.e., large $\rho$), this must be effectively counteracted by a strong enough damping component to ensure boundedness and global existence.
\end{remark}
\begin{remark}
Even though the model is naturally inspired by population dynamics and biological transport, its main interest lies in the mathematical challenges posed by the coexistence of local superlinear growth and nonlocal source terms, and the delicate balance needed to ensure solution boundedness and global well-posedness. 
\end{remark}
\section{Some preliminaries and auxiliary tools}\label{PreliminarySection}
In this section, we recall two fundamental and general results that serve as essential preliminaries for the proof of our main findings. These foundational statements provide the necessary theoretical framework and tools that underpin the subsequent analysis.

First, we will rely on the following variant of the Gagliardo--Nirenberg inequality.
\begin{lemma}\label{LemmaGN}
Let $\Omega$ satisfy the hypotheses in \eqref{reglocal}, and let $r\geq 1$  and $0<q\leq p \leq \infty$  fulfilling
$\frac{1}{r} \leq \frac{1}{n} + \frac{1}{p}. $
Then, for 
$\theta= \frac{\frac{1}{q}-\frac{1}{p}}{\frac{1}{q}+\frac{1}{n}-\frac{1}{r}},$
there exists $C_{GN}=C_{GN}(p,q,r,\Omega)>0$ such that for all $\psi \in W^{1,r}(\Omega) \cap L^q(\Omega)$ it holds that 
\[
\|\psi\|_{L^p(\Omega)} \leq C_{GN} (\|\nabla \psi\|_{L^r(\Omega)}^{\theta} \|\psi\|_{L^q(\Omega)}^{1-\theta}+ \|\psi\|_{L^q(\Omega)}).
\]
\begin{proof}
See \cite[Lemma 2.3]{LiLankeitNonLinearity}.
\end{proof}
\end{lemma}
We also make use of this result, consequence of Maximal Sobolev Regularity theory. 
\begin{lemma}\label{lem:MaxReg}
Let $\Omega$ satisfy the hypotheses in \eqref{reglocal}, and let  $q>1$. Then there is $C_P=C_P(q,\Omega)>0$ such that the following holds: Whenever $T\in(0,\infty]$, $I=[0,T)$, $h\in L^q(I;L^q(\Omega))$ and $\psi_0\in W^{2,q}_{\nu}= \{\psi_0\in W^{2,q}(\Omega)\,:\, \partial_\nu \psi_0=0 \,\textrm{ on }\, \partial\Omega\}$, every solution $\psi\in W_{loc}^{1,q}(I;L^q(\Omega))\cap L^q_{loc}(I;W^{2,q}(\Omega))$ of 
 \[
  \psi_t=\Delta \psi- \psi + h\;\; \text{ in }\;\;\Omega\times(0,T);\quad
  \partial_{\nu} \psi=0\;\; \text{ on }\;\;\partial\Omega \times(0,T); \quad \psi(\cdot,0)=\psi_0 \;\; \text{ on }\;\;\Omega  
 \]
 satisfies 
 \[
  \int_0^t  e^s{\int_{\Omega} |\Delta \psi(\cdot,s)|^q}ds \le C_P \left[\lVert \psi_0 \rVert_{q,1-\frac{1}{q}}^q+\int_0^t e^s {\int_{\Omega} |h(\cdot,s)|^{q}}ds\right] \quad \text{for all } t\in(0,T).
 \]
 \begin{proof}
The proof builds upon the classical result established in \cite{PrussSchanaubeltMaximalRegul}. For a suitable adaptation to the context considered here, see, for example, the detailed treatment in \cite{IshidaLankeitVigliloro-Gradient}.
 \end{proof}
\end{lemma}
\section{Local solutions and their main properties: a boundedness criterion}\label{SectionLocalInTime}
As it is well known in evolutionary models, the standard approach to establishing the boundedness of solutions involves deriving uniform-in-time a priori estimates for certain norms of the local solution to the problem under consideration. To this end, the first step is justifying the existence of local classical solutions to the model in examination.
\begin{lemma}[Local existence]\label{localSol}
For $\tau\in \{0,1\}$, let $\Omega$, $u_0, \tau v_0$ comply with hypotheses in \eqref{reglocal}.  Additionally, let $\chi, a, b, c>0$ and $\beta, \delta, \gamma, \rho \geq 1$.  
Then problem \eqref{problem} has a unique and nonnegative classical solution
\begin{equation*}
(u,v)\in C^{2+\alpha,1+\frac{\alpha}{2}}( \Bar{\Omega} \times [0, \TM))\times C^{2+\alpha,\tau+\frac{\alpha}{2}}( \Bar{\Omega} \times [0, \TM)),
\end{equation*}
for some maximal $\TM\in(0,\infty]$, obeying the following property:
 \begin{equation}\label{dictomyCriteC2+delVerions0}
 \text{either } \TM=\infty \quad \text{or}\quad \limsup_{t \to \TM} \left(\|u(\cdot,t)\|_{C^{2+\alpha}(\bar\Omega)}+\|v(\cdot,t)\|_{C^{2+\alpha}(\bar\Omega)}\right)=\infty.
\end{equation}
\begin{proof} 
We present a version of the existence proof which, despite keeping the rigorous structure, main technical tools, and reasoning, omits extended commentary; in \cite[Lemma 2.1]{IshidaLankeitVigliloro-Gradient} and \cite[Proposition 4]{BianEtAlNonlocal} details can be found.

For any $R>0$, we define the closed convex set
\[
S_T := \left\{ 0 \le u \in C^{1,\frac{\alpha}{2}}(\bar{\Omega} \times [0,T]) : \|u(\cdot,t) - u_0\|_{L^\infty(\Omega)} \le R,\ \forall t \in [0,T] \right\}.
\]
Given $\tilde{u} \in S_T$, the corresponding solution $v$ to the linear problem
\[
\begin{cases}
\tau v_t - \Delta v + v = \tilde{u} & \text{in } \Omega \times (0,T), \\
\partial_\nu v = 0 & \text{on } \partial \Omega \times (0,T), \\
\tau v(x,0) = \tau v_0(x) & x \in \bar{\Omega},
\end{cases}
\]
satisfies, do to bootstrap regularity arguments, $v \in C^{2+\alpha,\tau + \frac{\alpha}{2}}(\bar{\Omega} \times [0,T])$, with a norm bounded by a constant $H = H(R) > 0$ due to classical regularity theory.

Using with such a $v$ on our hands, we solve the nonlinear problem
\[
\begin{cases}
u_t = \Delta u - \chi \nabla \cdot (u \nabla v) + a u^\rho - c|\nabla u|^\gamma - b \left( \int_\Omega u^\beta \right)^\delta & \text{in } \Omega \times (0,T), \\
\partial_\nu u = 0 & \text{on } \partial \Omega \times (0,T), \\
u(x,0) = u_0(x) & x \in \bar{\Omega},
\end{cases}
\]
which admits a unique classical solution $u \in C^{2+\alpha,1+\frac{\alpha}{2}}(\bar{\Omega} \times [0,T])$ for some small $T>0$, due to Schauder theory. Moreover, one obtains
\[
\|u(\cdot,t) - u_0\|_{L^\infty(\Omega)} \le L(R) T^{1 + \frac{\alpha}{2}},
\]
so choosing $T \le \left(\frac{R}{L(R)}\right)^{\frac{2}{2+\alpha}}$ ensures $\Phi(\tilde{u}) := u \in S_T$. The map $\Phi: S_T \to S_T$ is compact by classical embedding results (cf.~\cite[Theorem 4.25]{BrezisBook}), and Schauder's fixed point theorem yields a fixed point $u \in S_T$.

By standard continuation arguments, this yields a maximal existence time $\TM \in (0,\infty]$ such that either the solution exists globally, or blow-up in $C^{2+\alpha,\tau+\frac{\alpha}{2}}$ occurs at $\TM$.
\end{proof}
\end{lemma}
From now on with $(u,v)$ we will refer to the unique local solution to model \eqref{problem} defined on $\Omega \times (0,\TM)$ and  provided by Lemma \ref{localSol}. 

The following result summarizes the parameter conditions under which our main claims (Theorem \ref{MainTheorem} and Corollary \ref{MainCor}) can be established. Indeed, the upcoming lemma ensures the boundedness of the mass, which, as thoroughly discussed in $\S$\ref{Intro}, is not automatic but rather depends on specific relations among the model parameters. This mass boundedness will later play a crucial role in controlling some $L^k(\Omega)$-norm of $u$ for large $k>1$.

(Let us remark that all constants \( c_i \), with \( i = 1, 2, \ldots \), appearing in the subsequent analysis are assumed to be positive.)
\begin{lemma}
Let $\rho> \beta\geq 1$, $\gamma >\frac{n\rho}{n+\beta}$ and $\delta>\max\left\{\frac{n\gamma-n\rho+\rho\gamma}{n\gamma-n\rho+\beta\gamma}, \frac{\rho}{\beta}\right\}$. Then there exists $M>0$ such that 
\begin{equation}\label{MassG}
\int_\Omega u(x,t)\,dx \leq M \quad \textrm{for all } t \in (0,\TM).
\end{equation}
If, indeed, $\delta=\frac{\rho}{\beta}$ and $\gamma=\rho\geq 1$, relation \eqref{MassG} is complied with $M=\int_\Omega u_0(x)dx$ whenever  
\begin{equation}\label{Cond_b}
\rho >\beta \quad \textrm{and} \quad b \geq \frac{a 2^{\rho-1}C_{GN}^\rho}{|\Omega|} \left[\frac{\beta \rho}{n\rho+\beta \rho-n \beta} \left(\frac{a 2^{\rho-1}C_{GN}^\rho n (\rho-\beta)}{c(n\rho+\beta \rho-n \beta)}\right)^{\frac{n(\rho-\beta)}{\beta \rho}}+1\right],
\end{equation}
being $C_{GN}=C_{GN}(\rho,\gamma,\beta,\Omega)$ the constant provided by Lemma \ref{LemmaGN}.
\begin{proof}
Integrating the first equation in \eqref{problem} over $\Omega$ yields
\begin{equation}\label{M01}
\frac{d}{dt} \int_\Omega u 
= a \int_\Omega u^\rho - c \int_\Omega |\nabla u|^\gamma 
- b |\Omega| \left(\int_\Omega u^\beta\right)^\delta 
\quad \text{for all } t \in (0,\TM).
\end{equation}
To estimate the term involving $\int_\Omega u^\rho$, we apply the Gagliardo--Nirenberg (Lemma \ref{LemmaGN}) and the Young inequalities. For $\rho > \beta$ and $\gamma > \frac{n\rho}{n + \beta}$, we obtain for some $\const{GNY}=\const{GNY}(C_{GN}), \const{c}=\const{c}(C_{GN})$ (we will omit to specify the dependence of other constants on $C_{GN}$ if not necessary)
\begin{equation}\label{GagliardoNirenbY}
\begin{split}
a \|u\|_{L^\rho(\Omega)}^\rho 
&\le \const{GNY} \|\nabla u\|_{L^\gamma(\Omega)}^{\rho \theta} \|u\|_{L^\beta(\Omega)}^{\rho(1-\theta)} + \const{GNY} \|u\|_{L^\beta(\Omega)}^\rho \\
&\le c \int_\Omega |\nabla u|^\gamma + \const{c} \left(\int_\Omega u^\beta\right)^{\frac{\rho \gamma (1 - \theta)}{\beta (\gamma - \theta \rho)}} + \const{GNY} \left(\int_\Omega u^\beta\right)^{\frac{\rho}{\beta}}\quad \text{for all } t \in (0,\TM),
\end{split}
\end{equation}
where 
\[
\theta:= \frac{-\frac{1}{\rho} + \frac{1}{\beta}}{\frac{1}{\beta} + \frac{1}{n} - \frac{1}{\gamma}} \in (0,1).
\]
Substituting \eqref{GagliardoNirenbY} into \eqref{M01} and applying Young’s inequality, we obtain, under the condition
\[
\delta > \max\left\{ \frac{n\gamma - n\rho + \rho \gamma}{n\gamma - n\rho + \beta \gamma}, \frac{\rho}{\beta} \right\},\]
the estimate
\begin{equation}\label{MM}
\frac{d}{dt} \int_\Omega u 
\le - \frac{b |\Omega|}{2} \left(\int_\Omega u^\beta\right)^\delta + \const{d} \quad \text{on } \;(0,\TM).
\end{equation}
Applying Hölder’s inequality for $\beta>1$ (the case $\beta=1$ does not require the estimate below), we further estimate
\[
\left(\int_\Omega u^\beta\right)^\delta \ge \const{G} \left( \int_\Omega u \right)^{\delta \beta},
\quad \text{with } t \in (0,\TM)\]
so that \eqref{MM} becomes
\[
\frac{d}{dt} \int_\Omega u \le \const{d} - \const{G} \left( \int_\Omega u \right)^{\delta \beta} 
\quad \text{for all } t \in (0,\TM),
\]
with initial condition
\[
\left(\int_\Omega u(x,t)\,dx\right)\Bigg|_{t=0} = \int_\Omega u_0(x)\,dx.
\]
An application of a standard comparison principle for ordinary differential inequalities then yields boundedness of the total mass on $(0,\TM)$.

Concerning the second assertion, we observe that for $\delta = \frac{\rho}{\beta}$ and $\gamma = \rho$, relation \eqref{M01} reduces to
\begin{equation}\label{M01B}
\frac{d}{dt} \int_\Omega u = a \int_\Omega u^\rho - c \int_\Omega |\nabla u|^\rho - b |\Omega| \left( \int_\Omega u^\beta \right)^{\frac{\rho}{\beta}} 
\quad \text{for all } t \in (0,\TM).
\end{equation}
Since $\rho > \beta$, we refine estimate \eqref{GagliardoNirenbY} by explicitly analyzing the dependence on constants. This leads for all $t\in (0,\TM)$ to
\begin{equation}\label{GagliardoNirenbYB}
\begin{split}
a \|u\|_{L^{\rho}(\Omega)}^{\rho} 
&\leq c \int_\Omega |\nabla u|^\rho \\ & 
+ a 2^{\rho - 1} C_{GN}^\rho \left[
\frac{\beta \rho}{n\rho + \beta \rho - n \beta} 
\left( \frac{a 2^{\rho - 1} C_{GN}^\rho n (\rho - \beta)}{c (n\rho + \beta \rho - n \beta)} \right)^{\frac{n(\rho - \beta)}{\beta \rho}} + 1
\right] 
\left( \int_\Omega u^\beta \right)^{\frac{\rho}{\beta}}.
\end{split}
\end{equation}
By plugging \eqref{GagliardoNirenbYB} into \eqref{M01B}, and under the second condition in \eqref{Cond_b},
we conclude that
\begin{equation*}
\frac{d}{dt} \int_\Omega u 
\le 0 \quad \text{for all } t \in (0,\TM).
\end{equation*}
Therefore, the total mass is nonincreasing in time, and we deduce the estimate
\[
\int_\Omega u(x,t)\,dx \le \int_\Omega u_0(x)\,dx 
\quad \text{for all } t \in (0,\TM).
\]
\end{proof}
\end{lemma}
\begin{remark}
Contrary to what was established in the previous lemma, in the case $\rho \leq \beta$ the uniform-in-time boundedness of the mass is obtained straightforwardly; indeed, for  $\beta>\rho>1$ and $\delta>1$, starting from \eqref{M01}, applications of Young's and Hölder's inequalities provide 
\begin{equation*}
\begin{split}
\frac{d}{dt} \int_\Omega u 
&\leq \int_\Omega u^\beta+\const{AA} 
- b |\Omega| \left(\int_\Omega u^\beta\right)^\delta \\ & \leq -\frac{b|\Omega|}{2} \left(\int_\Omega u^\beta\right)^\delta+\const{BB}\leq -\const{AB}\left(\int_\Omega u\right)^{\delta\beta}+\const{BB}
\quad \text{for all } t \in (0,\TM).
\end{split}
\end{equation*}
(For the case $\beta$ and/or $\delta=1$ the reasoning is even more direct.)

If $\rho=\beta\geq 1$, and $\delta>1$, the previous inequality is reorganized as 
\begin{equation*}
\begin{split}
\frac{d}{dt} \int_\Omega u 
\leq a \int_\Omega u^\rho
- b |\Omega| \left(\int_\Omega u^\rho\right)^\delta  \leq -\const{CC} \left(\int_\Omega u\right)^{\rho \delta}+\const{BBA}
\quad \text{on } \;\; (0,\TM).
\end{split}
\end{equation*}
(For the limit case $\rho=\beta$ and $\delta=1$, it is seen that $\int_\Omega u(x,t)dx$ is controlled whenever $a\leq b |\Omega|$.)

In the present study, we do not extend the investigation of the uniform-in-time boundedness of solution boundedness, which typically relies on the control of mass, as our primary focus (delineated in §\ref{MotivationSection}) concerns regimes where such control cannot be presumed a priori. Nevertheless, it should be emphasized that in cases where the total mass remains uniformly bounded, the criteria ensuring uniform boundedness of solutions are considerably less stringent than those articulated in Theorem~\ref{MainTheorem}.
\end{remark}
As shown in the previous Lemma \ref{localSol}, relation \eqref{dictomyCriteC2+delVerions0} holds, implying that the solution to problem \eqref{problem} blows up in finite time $\TM$ in the $C^{2+\alpha}(\Omega)$-norm. To ensure that blow-up also occurs in the $L^\infty(\Omega)$-norm at time $\TM$, further assumptions on the growth of the nonlinear term $|\nabla u|^\gamma$, with $\gamma \geq 1$, are necessary. Specifically, we rely on a classical regularity result by Lieberman, which guarantees the Hölder continuity of bounded solutions to quasilinear parabolic equations under gradient nonlinearities satisfying the Hadamard growth condition (see, e.g., \cite[Chapters IV--VI]{LiebermanBook}). Under such structural conditions, boundedness in $C^{2+\alpha}(\Omega)$ implies boundedness in $L^\infty(\Omega)$ up to $\TM$. 

We formalize this below.
\begin{lemma}[Boundedness criterion] \label{ExtensionLemma}
For $1\leq \gamma \leq 2$, let $u\in L^\infty((0,\TM);L^\infty(\Omega))$; then we have $\TM=\infty$ and in particular $u\in L^\infty((0,\infty);L^\infty(\Omega))$. 
\begin{proof}
The claim is proved by contradiction. Suppose that $\TM < \infty$. It suffices to show that the uniform-in-time boundedness of $u$ in $L^\infty(\Omega)$ for $t \in (0,\TM)$ implies the boundedness of 
\[
\sup_{t\in [0,\TM]} \left( \| u(\cdot,t) \|_{C^{2+\alpha}(\bar{\Omega})} + \| v(\cdot,t) \|_{C^{2+\alpha}(\bar{\Omega})} \right).
\]
This would contradict the blow-up alternative given by \eqref{dictomyCriteC2+delVerions0}, and thus complete the proof.

First, the boundedness of $\| v(\cdot,t) \|_{C^{1+\alpha}(\bar{\Omega})}$ on $[0,\TM]$ follows from the assumption $u \in L^\infty((0,\TM); L^\infty(\Omega))$ and the regularity results for linear parabolic equations. More precisely, applying either \cite[Theorem~1.2]{lieberman_paper} (when $\tau = 1$) or \cite[Theorem~I.19.1]{friedman_pde} (when $\tau = 0$), in conjunction with the Sobolev embedding $W^{2,k}(\Omega) \hookrightarrow C^{1+\alpha}(\bar{\Omega})$ for arbitrarily large $k$, ensures that $v \in C^{1+\alpha}(\bar{\Omega})$ uniformly in time.

Letting $X = (x,t)$ and defining the vector field and source term
\[
A(X,u,\nabla u) := \nabla u - \chi u \nabla v(X), \qquad B(X,u,\nabla u) := a u^\rho - c |\nabla u|^\gamma - b \left( \int_\Omega u^\beta \right)^\delta,
\]
we consider the reformulated system
\begin{equation*}\label{ProblemAuxiliary}
\begin{cases}
u_t = \nabla \cdot A(X,u,\nabla u) + B(X,u,\nabla u), & X \in \Omega \times (0,\TM),\\
A(X,u,\nabla u) \cdot \nu = 0, & X \in \partial\Omega \times (0,\TM),\\
u(x,0) = u_0(x), & x \in \bar{\Omega}.
\end{cases}
\end{equation*}
For $\gamma \in [1,2]$, the source term $B(X,u,\nabla u)$ satisfies the Hadamard structural conditions analogous to \cite[condition (1.2c), Theorem~1.2]{lieberman_paper}. Therefore, we can invoke \cite[Theorem~1.2]{lieberman_paper} to conclude that $\| u(\cdot,t) \|_{C^{1+\alpha}(\bar{\Omega})}$ remains bounded on $[0,\TM]$. Then, employing regularity arguments as mentioned in Lemma~\ref{localSol}, together with the regularity of the initial data $(u_0, \tau v_0)$, we obtain
\[
\sup_{t \in [0,\TM]} \left( \| u(\cdot,t) \|_{C^{2+\alpha}(\bar{\Omega})} + \| v(\cdot,t) \|_{C^{2+\alpha}(\bar{\Omega})} \right) < \infty,
\]
which contradicts \eqref{dictomyCriteC2+delVerions0} and thus the proof is completed.
\end{proof}
\end{lemma} 
\section{Uniform-in-time boundedness of $u$ in $L^k(\Omega)$. Proof of the results}\label{EstimatesAndProofSection}
We will see later that the uniform-in-time boundedness of \( u \) is implied whenever \( u \in L^\infty\big((0, T_{\mathrm{max}}); L^k(\Omega)\big) \) for sufficiently large \( k \); thereafter, let us proceed below to derive certain \textit{a priori} integral estimates.
\subsection{A priori boundedness of $u$ for model \eqref{problem}}
To further analyze the qualitative properties of solutions to the model, aimed at establishing the desired boundedness, we now present several key a priori estimates involving the function \( u \).
\begin{lemma}\label{LemmaPriori} 
Let $n \in \mathbb{N}$, $\beta, \gamma, \delta, \rho \geq 1$  satisfy 
\begin{equation}\label{LemmaMT1}
1 \leq \beta < \frac{2(n+1)}{n}, \; \max\{2, \beta\} < \rho < \frac{2(n+1)}{n}, \; \frac{n\rho}{n+1}<\gamma \leq 2 \quad \textrm{and} 
\quad \delta>\max\left\{\frac{n\gamma-n\rho+\rho\gamma}{n\gamma-n\rho+\beta\gamma}, \frac{\rho}{\beta}\right\},    
\end{equation} 
or
\begin{equation}\label{LemmaMT2}
1 \leq \beta < \rho \leq 2, \quad \frac{2n}{n+1}<\gamma \leq 2 \quad \textrm{and} 
\quad \delta>\max\left\{\frac{n\gamma-n\rho+\rho\gamma}{n\gamma-n\rho+\beta\gamma}, \frac{\rho}{\beta}\right\}. 
\end{equation} 
Then, there exists $k_0\geq1$ such that $u\in L^\infty((0,\TM);L^k(\Omega))$ for all $k>k_0$.
\begin{proof}
We begin by fixing \( k_0 = 1 \); this value may be tacitly increased as needed to ensure the applicability of the forthcoming estimates. 

In order to derive appropriate integral bounds, we test the first equation in \eqref{problem} with the function \( k u^{k-1} \), and integrate by parts. This procedure yields, for all \( k > k_0 \) and for every \( t \in (0, T_{\max}) \),
\begin{equation}\label{Stima0}
\begin{split}
y'(t)=y' := \frac{d}{dt} \int_\Omega u^k & = k \int_\Omega u^{k-1} u_t = -k(k-1)\int_\Omega u^{k-2} |\nabla u|^2 
+ k(k-1)\chi \int_\Omega u^{k-1} \nabla u \cdot \nabla v \\
&\quad + ka \int_\Omega u^{k+\rho-1} 
- kc \int_\Omega u^{k-1} |\nabla u|^\gamma 
- kb \int_\Omega u^{k-1} \left( \int_\Omega u^\beta \right)^\delta \\
&= -\frac{4(k-1)}{k} \int_\Omega \left| \nabla u^{\frac{k}{2}} \right|^2 
- \chi(k-1) \int_\Omega u^k \Delta v 
+ ka \int_\Omega u^{k+\rho-1} \\
&\quad - kc \left( \frac{k - 1 + \gamma}{\gamma} \right)^{-\gamma} \int_\Omega \left| \nabla u^{\frac{k - 1 + \gamma}{\gamma}} \right|^\gamma 
- kb \int_\Omega u^{k-1} \left( \int_\Omega u^\beta \right)^\delta.
\end{split}
\end{equation}
In the analysis that follows, we distinguish between the two cases corresponding to \( \tau = 0 \) and \( \tau = 1 \).
\subsubsection*{$\triangleright$ The case $\tau=0$}
Exploiting the second equation in \eqref{problem}, for every \( k > k_0 \) we deduce the following differential inequality:
\begin{equation}\label{Stima1}
\begin{split}
y' \leq{} & -\frac{4(k-1)}{k} \int_\Omega \left|\nabla u^{\frac{k}{2}}\right|^2 
+ \chi (k-1) \int_\Omega u^{k+1} 
+ ka \int_\Omega u^{k+\rho-1} \\
& - kc \left(\frac{k-1+\gamma}{\gamma}\right)^{-\gamma} \int_\Omega \left|\nabla u^{\frac{k-1+\gamma}{\gamma}}\right|^\gamma 
- kb \int_\Omega u^{k-1} \left(\int_\Omega u^\beta \right)^\delta,
\quad t \in (0, T_{\max}).
\end{split}
\end{equation}
We first assume condition \eqref{LemmaMT1}. Since \(\rho > 2\), an application of Young's inequality yields
\begin{equation}\label{Young1}
\chi (k-1) \int_\Omega u^{k+1} \leq \int_\Omega u^{k+\rho-1} + \const{Y1} \quad \text{for all } t \in (0, T_{\max}).
\end{equation}
To estimate the term \(\int_\Omega u^{k+\rho-1}\), we apply Lemma \ref{LemmaGN}. Given that \(\gamma > \frac{n \rho}{n+1}\) and taking into account the mass bound \eqref{MassG}, a combination of the Gagliardo--Nirenberg (again Lemma \ref{LemmaGN}) and the Young inequalities implies, for all \(t \in (0,T_{\max})\),
\begin{equation}\label{GNrho}
\begin{split}
(ka + 1) \int_\Omega u^{k+\rho-1} 
&= (ka+1) \left\| u^{\frac{k-1+\gamma}{\gamma}} \right\|_{L^{\frac{\gamma(k+\rho-1)}{k-1+\gamma}}(\Omega)}^{\frac{\gamma(k+\rho-1)}{k-1+\gamma}} \\
&\leq \const{GN32} \left\| \nabla u^{\frac{k-1+\gamma}{\gamma}} \right\|_{L^\gamma(\Omega)}^{\frac{\gamma(k+\rho-1)}{k-1+\gamma} \theta_1} 
\left\| u^{\frac{k-1+\gamma}{\gamma}} \right\|_{L^{\frac{\gamma}{k-1+\gamma}}(\Omega)}^{\frac{\gamma(k+\rho-1)}{k-1+\gamma} (1-\theta_1)} + \const{GN32} \left\| u^{\frac{k-1+\gamma}{\gamma}} \right\|_{L^{\frac{\gamma}{k-1+\gamma}}(\Omega)}^{\frac{\gamma(k+\rho-1)}{k-1+\gamma}} \\
&\leq \const{GN33} \left( \int_\Omega \left|\nabla u^{\frac{k-1+\gamma}{\gamma}} \right|^\gamma \right)^{\frac{(k+\rho-1) \theta_1}{k-1+\gamma}} + \const{GN331} \leq kc \left(\frac{k-1+\gamma}{\gamma}\right)^{-\gamma} \int_\Omega \left|\nabla u^{\frac{k-1+\gamma}{\gamma}} \right|^\gamma + \const{GN3},
\end{split}
\end{equation}
where 
\[
\theta_1 := \frac{\frac{k-1+\gamma}{\gamma} - \frac{k-1+\gamma}{\gamma(k+\rho-1)}}{\frac{k-1+\gamma}{\gamma} - \frac{1}{\gamma} + \frac{1}{n}} \in (0,1) \quad \text{for } k > k_0.
\]
Substituting estimates \eqref{Young1} and \eqref{GNrho} into \eqref{Stima1} yields
\begin{equation}\label{uk_estimate4}
y' \leq - \frac{4(k-1)}{k} \int_\Omega \left|\nabla u^{\frac{k}{2}}\right|^2 + \const{SI} \quad \text{for } t \in (0, T_{\max}).
\end{equation}
Alternatively, suppose now condition \eqref{LemmaMT2}. Since \(\rho \leq 2\), Young's inequality implies
\begin{equation}\label{Young2}
ka \int_\Omega u^{k+\rho-1} \leq \int_\Omega u^{k+1} + \const{Y2} \quad \text{for all } t \in (0, T_{\max}),
\end{equation}
where the inequality is not necessary if \(\rho = 2\).

Because \(\gamma > \frac{2n}{n+1}\), another application of the Gagliardo--Nirenberg and the Young inequalities, along with the boundedness of the mass \eqref{MassG}, implies on \((0, T_{\max})\),
\begin{equation}\label{GNk+1}
\begin{split}
(\chi(k-1) + 1) \int_\Omega u^{k+1} 
&= (\chi(k-1) + 1) \left\| u^{\frac{k-1+\gamma}{\gamma}} \right\|_{L^{\frac{\gamma(k+1)}{k-1+\gamma}}(\Omega)}^{\frac{\gamma(k+1)}{k-1+\gamma}} \\
&\leq \const{GN32a} \left\| \nabla u^{\frac{k-1+\gamma}{\gamma}} \right\|_{L^\gamma(\Omega)}^{\frac{\gamma(k+1)}{k-1+\gamma} \theta_2}
\left\| u^{\frac{k-1+\gamma}{\gamma}} \right\|_{L^{\frac{\gamma}{k-1+\gamma}}(\Omega)}^{\frac{\gamma(k+1)}{k-1+\gamma} (1-\theta_2)}  + \const{GN32a} \left\| u^{\frac{k-1+\gamma}{\gamma}} \right\|_{L^{\frac{\gamma}{k-1+\gamma}}(\Omega)}^{\frac{\gamma(k+1)}{k-1+\gamma}} \\
&\leq \const{GN33a} \left( \int_\Omega \left| \nabla u^{\frac{k-1+\gamma}{\gamma}} \right|^\gamma \right)^{\frac{(k+1) \theta_2}{k-1+\gamma}} + \const{GN33ac} + kc \left(\frac{k-1+\gamma}{\gamma}\right)^{-\gamma} \int_\Omega \left|\nabla u^{\frac{k-1+\gamma}{\gamma}} \right|^\gamma + \const{GN3a},
\end{split}
\end{equation}
where
\[
\theta_2 := \frac{\frac{k-1+\gamma}{\gamma} - \frac{k-1+\gamma}{\gamma(k+1)}}{\frac{k-1+\gamma}{\gamma} - \frac{1}{\gamma} + \frac{1}{n}} \in (0,1) \quad \text{for } k > k_0.
\]
Inserting inequalities \eqref{Young2} and \eqref{GNk+1} into \eqref{Stima1} recovers the estimate \eqref{uk_estimate4} up to a constant factor.

Finally, again the Gagliardo--Nirenberg and the Young inequalities, together with the mass bound \eqref{MassG}, imply for all \( t \in (0, T_{\max}) \),
\begin{equation}\label{GN2}
\begin{split}
\int_\Omega u^k 
&= \left\| u^{\frac{k}{2}} \right\|_{L^2(\Omega)}^2 
\leq \const{GN23} \left\| \nabla u^{\frac{k}{2}} \right\|_{L^2(\Omega)}^{2 \theta_3} \left\| u^{\frac{k}{2}} \right\|_{L^{\frac{2}{k}}(\Omega)}^{2 (1-\theta_3)} + \const{GN23} \left\| u^{\frac{k}{2}} \right\|_{L^{\frac{2}{k}}(\Omega)}^{2} \\
&\leq \const{GN24} \left( \int_\Omega \left| \nabla u^{\frac{k}{2}} \right|^2 \right)^{\theta_3} + \const{GN25} \leq \frac{4(k-1)}{k} \int_\Omega \left| \nabla u^{\frac{k}{2}} \right|^2 + \const{GN255},
\end{split}
\end{equation}
where 
\[
\theta_3 := \frac{\frac{k}{2} - \frac{1}{2}}{\frac{k}{2} - \frac{1}{2} + \frac{1}{n}} \in (0,1) \quad \text{for } k > k_0.
\]
Hence, combining \eqref{GN2} with \eqref{uk_estimate4} yields the differential inequality
\[
y' \leq \const{AC1} - y, \quad t \in (0, T_{\max}), \quad y(0)=y_0= \left(\int_\Omega u^k(x,t)\, dx \right){\Bigg|_{t=0}}=\int_\Omega u_0^k(x)\, dx,
\]
and by applying a standard comparison principle for ordinary differential equations, we conclude that
\[
\sup_{t \in (0, T_{\max})} \int_\Omega u^k \leq \max\{y_0,\const{AC1}\}.
\]
\subsubsection*{$\triangleright$ The case $\tau=1$}
Let us examine the integral term $-\chi (k-1) \int_\Omega u^k \Delta v$; it can be estimated by applying Young's inequality as
\begin{equation}\label{ParabSt1}
-\chi (k-1) \int_\Omega u^k \Delta v \leq \chi (k-1) \int_\Omega u^k |\Delta v| 
\leq \int_\Omega u^{k+1} + \const{GS} \int_\Omega |\Delta v|^{k+1} \quad \text{for all } t \in (0,T_{\max}).
\end{equation}
Incorporating estimate \eqref{ParabSt1} into \eqref{Stima0}, we can rewrite the latter as
\begin{equation}\label{ParabSt3}
\begin{split}
y' \leq\;& -\frac{4(k-1)}{k} \int_\Omega \left|\nabla u^{\frac{k}{2}}\right|^2 + \int_\Omega u^{k+1} + \const{GS} \int_\Omega |\Delta v|^{k+1} + k a \int_\Omega u^{k+\rho -1} \\
& - k c \left(\frac{k-1+\gamma}{\gamma}\right)^{-\gamma} \int_\Omega \left|\nabla u^{\frac{k-1+\gamma}{\gamma}}\right|^\gamma 
- k b \int_\Omega u^{k-1} \left(\int_\Omega u^\beta\right)^\delta
\quad \text{on } (0, T_{\max}).
\end{split}
\end{equation}
Next, we add the term $y=y(t)$ to both sides of \eqref{ParabSt3} and multiply through by $e^t$. Observing that
\[
e^t y' + e^t y = \frac{d}{dt} \left( e^t y \right),
\]
integration over the interval $(0,t)$ for any $t \in (0, T_{\max})$ leads to
\begin{equation}\label{StimaParab}
\begin{split}
e^t y \leq y(0) + \int_0^t e^s \Bigg( & -\frac{4(k-1)}{k} \int_\Omega \left|\nabla u^{\frac{k}{2}}\right|^2 + \int_\Omega u^{k+1} + \const{GS} \int_\Omega |\Delta v|^{k+1} + k a \int_\Omega u^{k+\rho -1} \\
& - k c \left(\frac{k-1+\gamma}{\gamma}\right)^{-\gamma} \int_\Omega \left|\nabla u^{\frac{k-1+\gamma}{\gamma}}\right|^\gamma
- k b \int_\Omega u^{k-1} \left(\int_\Omega u^\beta \right)^\delta \Bigg) ds.
\end{split}
\end{equation}
At this stage, by applying Lemma~\ref{lem:MaxReg} to the second equation of \eqref{problem} with $\psi = v$, $q = k+1$, and $h = u$, we obtain the estimate
\begin{equation}\label{ParabSt2}
\const{GS} \int_0^t e^s \left( \int_\Omega |\Delta v(\cdot,s)|^{k+1} \right) ds \leq \const{G_2} + \const{G_3} \int_0^t e^s \left( \int_\Omega u^{k+1} \right) ds \quad \text{for all } t \in (0, T_{\max}),
\end{equation}
where $\const{G_2}, \const{G_3}$ depend also on $C_P$. (We will not explicitly state the dependence of other constants on $C_P$, unless required.) 

Substituting \eqref{ParabSt2} into \eqref{StimaParab} and proceeding analogously to the case $\tau=0$, particularly in handling the terms $\int_\Omega u^{k+1}$ and $\int_\Omega u^{k+\rho -1}$, we arrive at
\[
e^t \int_\Omega u^k \leq \const{sd} + \const{a_6} (e^t - 1) \quad \text{for all } t \in (0, T_{\max}),
\]
thus concluding the proof, having naturally $\int_\Omega u^k \leq \const{sd1}$ on $(0,\TM).$
\end{proof}
\end{lemma}
\begin{lemma}\label{LemmaPrioriCor}
Let 
$a, c >0$, $\beta, \gamma=\rho \geq 1, \delta=\frac{\rho}{\beta}$ fulfill   
\begin{equation*}\label{MCC}
\max\left\{\beta, \frac{2n}{n+1}\right\} < \rho \leq 2 \quad \textrm{and} \quad  b \geq \frac{a 2^{\rho-1}C_{GN}^\rho}{|\Omega|} \left[\frac{\beta \rho}{n\rho+\beta \rho-n \beta} \left(\frac{a 2^{\rho-1}C_{GN}^\rho n (\rho-\beta)}{c(n\rho+\beta \rho-n \beta)}\right)^{\frac{n(\rho-\beta)}{\beta \rho}}+1\right],
\end{equation*}
where $C_{GN}=C_{GN}(\rho,\gamma,\beta,\Omega)$ is the positive constant given in Lemma \ref{LemmaGN}. Then, there exists $k_1\geq1$ such that $u\in L^\infty((0,\TM);L^k(\Omega))$ for all $k>k_1$. 
\begin{proof}
Let us consider $k_1=1$; as done before, we will enlarge this initial value when necessary.
By taking $\gamma=\rho$ and $\delta=\frac{\rho}{\beta}$, relation \eqref{Stima0} becomes for all $k>k_1$
\begin{equation*}\label{Stima0C}
\begin{split}
y' =&-\frac{4(k-1)}{k} \int_\Omega |\nabla u^{\frac{k}{2}}|^2 -\chi (k-1) \int_\Omega u^k \Delta v+ka \int_\Omega u^{k+\rho-1}\\
&-kc \left(\frac{k-1+\rho}{\rho}\right)^{-\rho} \int_\Omega |\nabla u^{\frac{k-1+\rho}{\rho}}|^\rho - kb \int_\Omega u^{k-1} \left(\int_\Omega u^{\beta}\right)^{\frac{\rho}{\beta}}
\quad  \mbox{for all }\ t \in (0, \TM). 
\end{split}
\end{equation*}
We distinguish the cases $\tau=0$ and $\tau=1$.
\subsubsection*{$\triangleright$ The case $\tau=0$}
By exploiting the second equation in \eqref{problem} and relying on Young's inequality in view of $\rho<2$, we have on $(0, \TM)$ for all $k>k_1$
\begin{equation}\label{Stima1C}
y' \leq -\frac{4(k-1)}{k} \int_\Omega |\nabla u^{\frac{k}{2}}|^2 +(\chi (k-1)+1) \int_\Omega u^{k+1} -kc \left(\frac{k-1+\rho}{\rho}\right)^{-\rho} \int_\Omega |\nabla u^{\frac{k-1+\rho}{\rho}}|^\rho + \const{CC1}.
\end{equation}
(If $\rho=2$ this is immediate).
Similarly to what have done in the proof of Lemma \ref{LemmaPriori}, we arrive at
\begin{equation*}
y' \leq - \frac{4(k-1)}{k}  \int_\Omega |\nabla u^{\frac{k}{2}}|^2 + \const{siC}  \quad \mbox{on}\ (0, \TM).
\end{equation*}
At this stage the conclusion follows by reasoning as in the proof of Lemma \ref{LemmaPriori}.
\subsubsection*{$\triangleright$ The case $\tau=1$}
By recalling the values of $\gamma$ and $\delta$, estimate \eqref{StimaParab} can be rewritten for all $t \in (0,\TM)$ as
\begin{equation*}\label{StimaParabC}
\begin{split}
e^t y \leq y(0) +\int_0^t e^s &\left(-\frac{4(k-1)}{k} \int_\Omega |\nabla u^{\frac{k}{2}}|^2 + \int_\Omega u^{k+1} +  \const{GS} \int_\Omega |\Delta v|^{k+1}
 +ka \int_\Omega u^{k+\rho-1} \right.\\
& \quad \left. -kc \left(\frac{k-1+\rho}{\rho}\right)^{-\rho} \int_\Omega |\nabla u^{\frac{k-1+\rho}{\rho}}|^\rho - kb \int_\Omega u^{k-1} \left(\int_\Omega u^{\beta}\right)^{\frac{\rho}{\beta}}\right)\, ds.
\end{split}
\end{equation*}
By inserting in the above estimate bound \eqref{ParabSt2} and following the argument in the proof of Lemma \ref{LemmaPriori}, we obtain the claim.
\end{proof}
\end{lemma}

We can now conclude and give the proofs of the claimed results.
\subsection{Proofs of Theorem \ref{MainTheorem} and Corollary \ref{MainCor}}\label{SectionProofTheorem}
By virtue of assumptions \eqref{MT1} or \eqref{MT2} of Theorem \ref{MainTheorem}, or \eqref{MC} of Corollary \ref{MainCor}, Lemma \ref{LemmaPriori} and Lemma \ref{LemmaPrioriCor} hold for every $k>\bar{k}:=\max\{k_0, k_1\}$. 
Thereafter, $u \in L^{\infty}((0,\TM);L^k(\Omega))$, for every arbitrarily large $k$, and subsequently 
elliptic and parabolic regularity results, in conjunction with Sobolev embeddings, imply $\nabla v\in L^\infty((0,\TM);L^\infty(\Omega))$. In addition, the source $B(X,u,\nabla u)=au(x,t)^\rho-b\left(\int_\Omega  u(x,t)^\beta dx\right)^\delta-c|\nabla u(x,t)|^\gamma$ in problem \eqref{problem} is such that $B(X,u,\nabla u)\leq au(x,t)^\rho$ for all $(x,t)\in \Omega \times (0,\TM)$, so that $B(X,u,\nabla u)\in L^\infty((0,\TM);L^{k}(\Omega))$. 
As a consequence, from \cite[Lemma A.1]{TaoWinkParaPara} we have 
$u \in L^{\infty}((0,\TM);L^{\infty}(\Omega))$, so concluding by applying Lemma \ref{ExtensionLemma}. 
\qed

\subsubsection*{\bf\textit{\quad Acknowledgments}}
TL is supported by NNSF of P. R. China (Grant No. 61503171), CPSF (Grant No. 2015M582091), and NSF of Shandong Province (Grant No. ZR2016JL021).
SF and GV are members of the Gruppo Nazionale per l'Analisi Matematica, la Probabilit\`a e le loro Applicazioni (GNAMPA) of the Istituto Nazionale di Alta Matematica (INdAM) and they are partially 
supported by MIUR (Italian Ministry of Education, University and Research) Prin 2022 \emph{Nonlinear differential problems with applications to real phenomena} (Grant Number: 2022ZXZTN2). GV
is also partially supported by the research project {\em Partial Differential Equations and their role in understanding natural phenomena}, CUP F23C25000080007, funded by Fondazione di Sardegna annuality (2023) and \emph{Studio di modelli
nelle scienze della vita} (CUP J55F21004240001, ``DM 737/2021 risorse 2022–2023'', founded by European
Union–NextGenerationEU).


\end{document}